\documentclass[12pt,letterpaper]{article}
\usepackage{mathrsfs}
\usepackage{indentfirst}
\usepackage{makeidx}
\usepackage{amsmath,amssymb,amsfonts,amsthm,graphics,graphicx}
\usepackage{makeidx}
\usepackage{float}
\usepackage{color}

\newtheorem{thm}{Theorem}[section]

\newtheorem{cor}[thm]{Corollary}

\theoremstyle{definition}

\newtheorem{pro}[thm]{Proposition}

\def\-{\mbox{--}}

\newtheorem{remark}{Remark}

\def\pf{\noindent {\it Proof.} }

\begin{document}
\title{Distance proper connection of graphs}
\author{
\small Xueliang Li$^1$, Colton Magnant$^2$, Meiqin Wei$^1$, Xiaoyu Zhu$^1$\\
\small $^1$Center for Combinatorics and LPMC\\
\small Nankai University, Tianjin 300071, China\\
\small Email: lxl@nankai.edu.cn; weimeiqin8912@163.com; zhuxy@mail.nankai.edu.cn\\
\small $^2$Department of Mathematical Sciences\\
\small Georgia Southern University, Statesboro, GA 30460-8093, USA\\
\small Email: cmagnant@georgiasouthern.edu}
\date{}
\maketitle

\begin{abstract}
Let $G$ be an edge-colored connected graph. A path $P$ in $G$ is called a distance $\ell$-proper path if no two edges of the same color appear with fewer than $\ell$ edges in between on $P$. The graph $G$ is called $(k,\ell)$-proper connected if every pair of distinct vertices of $G$ are connected by $k$ pairwise internally vertex-disjoint distance $\ell$-proper paths in $G$. For a $k$-connected graph $G$, the minimum number of colors needed to make $G$ $(k,\ell)$-proper connected is called the $(k,\ell)$-proper connection number of $G$ and denoted by $pc_{k,\ell}(G)$. In this paper, we prove that $pc_{1,2}(G)\leq 5$ for any $2$-connected graph $G$. Considering graph operations, we find that $3$ is a sharp upper bound for the $(1,2)$-proper connection number of the join and the Cartesian product of almost all graphs. In addition, we find some basic properties of the $(k,\ell)$-proper connection number and determine the values of $pc_{1,\ell}(G)$ where $G$ is a traceable graph, a tree, a complete bipartite graph, a complete multipartite graph, a wheel, a cube or a permutation graph of a nontrivial traceable graph.

{\flushleft\bf Keywords}: distance $\ell$-proper path; $(k,\ell)$-proper connected; $(k,\ell)$-proper connection number

{\flushleft\bf AMS classification 2010}: 05C15, 05C40
\end{abstract}

\section{Introduction}
All graphs in this paper are finite, undirected, simple and connected. We follow the notation and terminology in the book \cite{BM}.

Based on the communication of information between agencies of the government, an immediate question is put forward as follows. What is the minimum number of passwords or firewalls needed that allows one or more secure paths between every two agencies so that the passwords along each path are distinct? This situation can be represented by a graph and studied by means of what is called rainbow colorings introduced by Chartrand et al.~in \cite{CJMZ}. An \emph{edge-coloring} of a graph is a mapping from its edge set to the set of natural numbers. A path in an edge-colored graph with no two edges sharing the same color is called a \emph{rainbow path}. A graph $G$ with an edge-coloring $c$ is said to be \emph{rainbow connected} if every pair of distinct vertices of $G$ is connected by at least one rainbow path in $G$. The coloring $c$ is called a \emph{rainbow coloring} of the graph $G$. For a connected graph $G$, the minimum number of colors needed to make $G$ rainbow connected is the rainbow connection number of $G$ and denoted by $rc(G)$. Many researchers have studied problems on the rainbow connection and got plenty of nice results, see \cite{KY,LS,LS2} for examples. For more details we refer to the survey paper \cite{LSS} and the book \cite{LS}.

Related to the rainbow connection number, researchers have become interested in the following question. What is the minimum number of passwords or firewalls that allows one or more secure paths between every two agencies where, as we progress from one agency to another along such a path, we are required to change passwords? Inspired by this, Borozan et al.~in \cite{BFG} and Andrews et al.~in \cite{ALLZ} introduced the concept of proper-path coloring of graphs. Let $G$ be an edge-colored graph. A path $P$ in $G$ is called a \emph{proper path} if no two adjacent edges of $P$ are colored with the same color. An edge-colored graph $G$ is \emph{$k$-proper connected} if every pair of distinct vertices $u,v$ of $G$ are connected by $k$ pairwise internally vertex-disjoint proper $(u,v)$-paths in $G$. For a connected graph $G$, the minimum number of colors needed to make $G$ $k$-proper connected is called the \emph{$k$-proper connection number} of $G$ and denoted by $pc_k(G)$. Particularly for $k=1$, we write $pc_1(G)$, the proper connection number of $G$, as $pc(G)$ for simplicity. Recently, many results have been obtained with respect to several aspects, such as connectivity, minimum degree, complements, operations on graphs and so on. For more details, we refer to the dynamic survey paper \cite{LM}.

We now describe a new parameter which serves as a bridge between proper connection and rainbow connection. A path $P$ in $G$ is called a \emph{distance $\ell$-proper path} if no two edges of the same color can appear with fewer than $\ell$ edges in between on $P$. The graph $G$ is called \emph{$(k,\ell)$-proper connected} if there is an edge-coloring $c$ such that every pair of distinct vertices of $G$ is connected by $k$ pairwise internally vertex-disjoint distance $\ell$-proper paths in $G$. Such a coloring is called a \emph{$(k,\ell)$-proper-path coloring} of $G$. In addition, if $t$ colors are used, then $c$ is referred to as a \emph{$(k,\ell)$-proper-path $t$-coloring} of $G$. For a connected graph $G$, the minimum number of colors needed to make $G$ $(k,\ell)$-proper connected is called the \emph{$(k,\ell)$-proper connection number} of $G$ and denoted by $pc_{k,\ell}(G)$.

From the definition, we can easily see that the $(k,1)$-proper connection number of a graph $G$ is actually its $k$-proper connection number, i.e., $pc_{k,1}(G)=pc_k(G)$. Meanwhile, the $(1,\ell)$-proper connection number of a graph $G$ can be its rainbow connection number as long as $\ell$ is large enough. Now we return to the example of communication of information between agencies of the government. The higher the required security level, the more the passwords or firewalls are needed, and the larger the cost. Interestingly, different security levels can be reflected by the values of $\ell$ in the $(k,\ell)$-proper connection number of a graph. Additionally, for $k=1$ and $\ell=2$, there is an edge-coloring using $pc_{1,2}$ colors such that there exists a $2$-proper path between each pair of vertices of the graph $G$. Furthermore, if we require that every path in $G$ is a $2$-proper path, then the edge-coloring becomes a strong edge-coloring. The strong chromatic index $\chi'_s(G)$, which was introduced by Fouquet and Jolivet \cite{FJ}, is the minimum number of colors needed in a  strong edge-coloring of $G$. Immediately we get that $pc_{1,2}(G)\leq \chi'_s(G)$. And this inspires us to pay particular attention to the $(1,2)$-proper connection number of a connected graph $G$, i.e., $pc_{1,2}(G)$.

This paper is organized as follows. In Section~\ref{Sect:Basic}, we find some basic properties of the $(k,\ell)$-proper connection number of a graph $G$ and determine the values of $pc_{1,\ell}(G)$ when $G$ is a traceable graph, a tree, a complete bipartite graph, a complete multipartite graph, a wheel or a cube. In Section~\ref{Sect:2-conn}, we study the $(1,2)$-proper connection number of $G$ for any $2$-connected graph $G$. Finally in Section~\ref{Sect:Oper}, we investigate the $(1,2)$-proper connection numbers of the join and the Cartesian product of graphs. Also, we obtain the value of $pc_{1,\ell}(G)$ in which $G$ is a permutation graph of a nontrivial traceable graph.

\section{Preliminaries}\label{Sect:Basic}

In this section, we first introduce some definitions and present several basic propositions, which will serve as essential tools in the following proof. Then the $(k,\ell)$-proper connection numbers of a series of simple graphs are characterized.

Before all of these, we introduce some basic symbols that will appear in the sequel. Let $G$ be a connected graph. We denote by $n$ the number of its vertices and $m$ the number of its edges. The \emph{distance between two vertices} $u$ and $v$ in $G$, denoted by $d(u,v)$, is the length of a shortest path between them in $G$. The \emph{eccentricity} of a vertex $v$ is $ecc(v):=max_{x\in V(G)}d(v, x)$. The \emph{radius} of $G$ is $rad(G):=min_{x\in V(G)}ecc(x)$. We write $\sigma'_2(G)$ as the largest sum of degrees of vertices $x$ and $y$, where $x$ and $y$ are taken over all couples of adjacent vertices in $G$. Additionally, we set $[n]=\{1,2,\dots,n\}$ for any integer $n\geq1$.

Similar to the case of $pc(G)=pc_{1,1}(G)$, we know that $pc_{1,\ell}(G)=1$ if and only if $G=K_{n}$. It can be easily checked that $pc_{1,\ell}(K_{n})=1$. The converse also holds since a noncomplete graph must have diameter greater than $1$, and as a result, at least two colors are needed. Apart from this, another two essential properties are posed as follows.

\begin{pro}\label{sr}
Let $G$ be a nontrivial connected graph and $\ell\geq 1$ be an integer, then we have $pc(G)\leq pc_{1,2}(G)\leq \cdots \leq pc_{1,\ell-1}(G)\leq pc_{1,\ell}(G)\leq rc(G)\leq m$.
\end{pro}
\pf From the definition of $(k,\ell)$-proper connection number, we know that $pc(G)=pc_{1,1}(G)$. In addition, $pc_{1,i-1}(G)\leq pc_{1,i}(G)$ for all $2\leq i\leq \ell$ since every distance $\ell$-proper path is also an distance $(\ell-1)$-proper path. What's more, a rainbow path is no doubt a $t$-proper path for all $t$ with $1\leq t\leq \ell$ and $rc(G)$ must be no greater than $m$.\qed

\begin{pro}\label{spanning}
If $G$ is a nontrivial connected graph, $H$ is a connected spanning subgraph of $G$ and $\ell\geq 1$ is an integer, then $pc_{1,\ell}(G)\leq pc_{1,\ell}(H)$. Particularly, $pc_{1,\ell}(G)\leq pc_{1,\ell}(T)$ for every spanning tree $T$ of $G$.
\end{pro}
\pf Certainly a $(1,\ell)$-proper-path coloring of $H$ can be extended to a $(1,\ell)$-proper-path coloring of $G$ by arbitrarily assigning used colors to the remaining edges. This implies $pc_{1,\ell}(G)\leq pc_{1,\ell}(H)$. \qed

We also give the following result for the $(1,\ell)$-proper connection number of the traceable graphs, i.e., graphs containing a hamiltonian path.

\begin{pro}\label{traceable}
Let $G$ be a traceable graph and $\ell$ be a positive integer, then $pc_{1,\ell}(G)\leq\ell+1$. Particularly, $pc_{1,2}(G)\leq3$.
\end{pro}
\pf  We simply consider a hamiltonian path of $G$ and color the edges in the sequence of $1,2,\dots,\ell,\ell+1,1,2,\dots,\ell,\ell+1,\dots$. \qed

Next we discuss the $(1,\ell)$-proper connection numbers of the trees, complete bipartite graphs, complete multipartite graphs, wheels and cubes. First, we show that the $(1,2)$-proper connection number of a tree $T$ is closely related to $\sigma'_{2}(T)$.
\begin{thm}\label{tree}
If $T$ is a nontrivial tree, then $pc_{1,2}(T)=\sigma'_2(T)-1$.
\end{thm}
\pf It is known that in a tree, every path is the unique path connecting its two end vertices, so we have to make sure that every path in this tree is a distance $2$-proper path. Let $x$ and $y$ be two vertices such that $d(x)+d(y)=\sigma'_2(T)$. Then we define the edge $xy$ to be \emph{level} $0$ and give an edge $e$ \emph{level} $i$ if any path containing $e$ and $xy$ has exactly $i-1$ edges in between. We also define the \emph{parent} of each vertex $v$ to be the vertex adjacent to it on the path from $v$ to $x$(or $y$). Thus every vertex $v$ except $x$ and $y$ has a unique parent, denoted by $p(v)$. A \emph{child} of a vertex $v$, denoted by $c(v)$, is a vertex $w$ such that $p(w) = v$.

Next we begin to give our coloring. The colors of edges incident to either $x$ or $y$ must be pairwise different because every pair of these edges lies on a path with fewer than $2$ edges in between. This means that $T$ requires at least $\sigma'_{2}(T)-1$ colors and so $pc_{1,2}(T)\geq \sigma'_2(T)-1$. Note that this assigns distinct colors to all edges of levels $0$ and $1$. Suppose we have given colors to all edges of levels $k\leq i \ (i\geq 1)$, now randomly choose an edge $e=wc(w)$ of level $i+1$. Obviously all edges incident to $c(w)$ have not been given colors yet since they all have level $i+2$ except for $e$. Thus we only need to consider the colored edges which are incident to either $w$ or $p(w)$. And the total number of all these edges is at most $\sigma'_{2}(T)-1$. In this way we can choose a color for $e$ that does not appear on an edge incident to $w$ or $p(w)$. Consequently, we can color all edges of level $i+1$. Repeating this procedure, the coloring of $T$ is well defined. One can easily check this is a $2$-proper-path coloring of $T$, which proves that $pc_{1,2}(T)=\sigma'_2(T)-1$. \qed

For a nontrivial tree $T$, the proof of the above theorem suggests the method to find the value of $pc_{1,\ell}(T)$ for general $\ell$. From the set of subtrees with diameter $\ell+1$, we choose a subtree $T_0$ with maximum size and color all its edges with distinct colors. If $\ell$ is odd, then we take all edges in $T_0$ adjacent to some $v_1\in T_0$ as \emph{level $0$} of $T$, where $ecc_{T_0}(v_1)=\frac{\ell+1}{2}$. If $\ell$ is even, then we choose some edge $u_0v_0\in T_0$ as \emph{level $0$} of $T$, where $ecc_{T_0}(u_0)=ecc_{T_0}(v_0)=\frac{\ell}{2}+1$. Similarly, we divide $V(T)\backslash V(T_0)$ into different levels and color the edges of $T$ successively. As a result, we can obtain the following corollary.
\begin{cor}\label{corollary1}
For $T$ a nontrivial tree, $pc_{1,\ell}(T)$ is equal to the maximum size of subtrees with diameter $\ell+1$ of $T$.
\end{cor}

For $\ell=1$, we have $pc_{1,1}(T)=pc(T)=\Delta(T)$ for any tree $T$. Together with Theorem \ref{tree}, we get the following corollary.

\begin{cor}\label{corollary2}
For every pair of integers $a,b$ where $1\leq a\leq b\leq 2a-1$, there exists a connected graph $G$ such that $pc(G)=a$ and $pc_{1,2}(G)=b$.
\end{cor}

We now turn our attention to the complete bipartite graph $K_{m,n}$. Let $K_{m, n} = U \cup V$ where $|U|=m$ and $|V|=n$. Without loss of generality, we assume that $m\leq n$.
\begin{thm}\label{bipartite}
Let $\ell \geq 2$ be an integer and $m\leq n$. Then,
\begin{eqnarray*}
pc_{1,\ell}(K_{m,n})=\left\{
\begin{array}{rcl}
n     &      & if\ m=1,\\
2     &      & if\ m\geq2\ and\ m\leq n\leq2^m,\\
3     &      & if\ \ell=2,\ m\geq2\ and\ n>2^m,\ or\\
~     &      & \ell\geq3,\ m\geq2\ and\ 2^m<n\leq3^m,\\
4     &      & if\ \ell\geq3,\ m\geq2\ and\ n>3^m.
\end{array} \right.
\end{eqnarray*}
\end{thm}
\pf Set $U=\{u_1,\dots,u_m\}$ and $|V|=n$.

If $m=1$, we can easily get that any two edges of $K_{1,n}$ must be assigned different colors in order to make $K_{1,n}$ proper connected, thus $pc_{1,\ell}(K_{1,n})=n$.

If $m\geq2$ and $m\leq n\leq2^m$, we define a $(1,\ell)$-proper-path $2$-coloring of $K_{m,n}$ for any $\ell\geq 2$ as follows. For each element $v\in V$, we assign a vector $v'=(v'_1,v'_2,\dots,v'_m)$ to it such that $v'_i\in \{1,2\},~i=1,\dots,m$ and vectors $(2,1,\dots,1)$, $(1,2,\dots,1),\cdots,(1,\dots,1,2)$ are all present. In addition, we make sure that $v'\neq w'$ for any distinct vertices $v,w\in V$. Color the edge $vu_i$ with $v'_i$. Thus for any $x,y\in V$, there exists $1\leq i\leq m$ such that $x'_i\neq y'_i$, so $xu_iy$ is a distance $\ell$-proper $(x,y)$-path. And for $u_i,u_j\in U$, $u_iwu_j$ is a distance $\ell$-proper $(u_i,u_j)$-path, where $w\in V$ and $w'=(w'_1,\dots,w'_{i-1},w'_i,w'_{i+1},\dots,w'_m)
=(1,\dots,1,2,1,\dots,1)$. For $u_t\in U$ and $z\in V$, the edge $u_tz$ is trivially a distance $\ell$-proper path in $K_{m,n}$. Hence, $pc_{1,\ell}(K_{m,n})=2$ for all $\ell\geq 2$.

If $m\geq2$, $n>2^m$ and $\ell=2$, we define a $(1,2)$-proper-path $3$-coloring of $K_{m,n}$ as follows. Suppose $V=V_1\cup V_2$ such that $|V_1|=2^m$. Similar to the previous case, for any $v\in V_1$, we assign a vector $v'=(v'_1,v'_2,\dots,v'_m)$ to $v$ such that $v'_i\in \{1,2\},~i=1,\dots,m$ and $v'\neq w'$ for any distinct vertices $v,w\in V_1$. In addition, for any vertex $t\in V_2$, we define $t'=(t'_1,\dots,t'_m)=(3,\dots,3)$ as its corresponding vector. Then for any $x\in V$, color the edge $xu_i$ with $x'_i$. Similar to the case above, for any pair of vertices in $V_{1} \cup U$, there is a distance $2$-proper path in between. For all $x,y\in V_2$, we see that $xu_1zu_2y$ is a distance $2$-proper $(x,y)$-path, where $z\in V_1$ and $z'=(1,2,1,\dots,1)$. And for any $x\in V_1$, $y\in V_2$, $xu_1y$ is a distance $2$-proper $(x,y)$-path. Thus $pc_{1,2}(K_{m,n})\leq 3$. For sharpness, we need to show that $pc_{1,2}(K_{m,n})>2$. If not, there is a $(1,2)$-proper $2$-coloring $c$ of $K_{m,n}$. However, according to the pigeon hole principle, there must exist $x,y\in V$ such that $c(xu_i)=c(yu_i)$ for $1\leq i\leq m$. Thus $y$ can not be reached from $x$ through a distance $2$-proper path of length two. Since any distance $2$-proper path of length at least $3$ requires the use of at least three colors, this contradicts the assumption that $pc_{1,2}(K_{m,n})=2$. Therefore, we have $pc_{1,2}(K_{m,n})=3$.

If $m\geq2$, $2^m<n\leq3^m$ and $\ell\geq3$, since $n>2^m$, again we know that two colors are not enough. We then define a $(1,\ell)$-proper-path $3$-coloring as follows. As above, we assign to the vertices of $V$ distinct vectors of length $m$ with entries from $\{1,2,3\}$ such that the vectors $(2,1,\dots,1)$, $(1,2,\dots,1),\cdots,(1,\dots,1,2)$ are all present. Obviously this is a $(1,\ell)$-proper-path coloring of $K_{m,n}$, which implies that $pc_{1,\ell}(K_{m,n})=3$.

Finally suppose $m\geq2$, $n>3^m$ and $\ell\geq3$. Since $n>3^m$, by the pigeonhole principle, there must be two vertices $x$ and $y$ in $U$ such that the edges $wx$ and $wy$ have the same color for all $w \in V$. Any $\ell$-proper path from $x$ to $y$ must have length at least $4$ so with $\ell \geq 3$, we know that three colors are not enough, so $pc_{1, \ell}(K_{m, n}) \geq 4$. We then define a $(1,\ell)$-proper-path $4$-coloring as follows. Suppose $V=V_1\cup V_2$ such that $|V_1|=2^m$. Similarly, we assign to the vertices of $V_{1}$ distinct vectors of length $m$ with entries from $\{1,2\}$. For any vertex $w\in V_{2}$, we assign to it the corresponding vector $w'=(3,4,4,\dots,4,4)$. With similar analysis, we know that for any $x,y\in V_1$ and $u_i,u_j\in U$, they both have distance $\ell$-proper paths in between. For any $x,y\in V_2$, $xu_1zu_2y$ is a distance $\ell$-proper $(x,y)$-path, where $z\in V_1$ and $z'=(1,2,1,\dots,1)$. And for any $x\in V_1$, $y\in V_2$, $xu_1y$ is a distance $\ell$-proper $(x,y)$-path. Thus $pc_{1,\ell}(K_{m,n})=4$. And the proof is completed.\qed

Next we consider a more general problem of determining $pc_{1,\ell}(G)$ when $G$ is a complete multipartite graph. For $t\geq3$, let $K_{n_1,\cdots,n_t}$ be the complete multipartite graph with $1\leq n_1\leq\cdots\leq n_t$. In addition, we set $m=\sum_{i=1}^{t-1}n_i$ and $n=n_t$. Then we have the following conclusion.

\begin{thm}\label{multipartite}
Let $t\geq3$, $\ell \geq 1$, $1\leq n_1\leq\cdots\leq n_t$, $m=\sum_{i=1}^{t-1}n_i$ and $n=n_t$. Then,
\begin{eqnarray*}
pc_{1,\ell}(K_{n_1,\cdots,n_t})=\left\{
\begin{array}{rcl}
1     &      & if\ n=1,\\
2     &      & if\ 2\leq n\leq2^m,\\
3     &      & otherwise.
\end{array} \right.
\end{eqnarray*}
\end{thm}
\pf For convenience, we write $G$ for $K_{n_1,\cdots,n_t}$. Let $V_1,\cdots,V_t$ be the partition of $G$ where $|V_{i}|=n_{i}$ for all $i$ with $1 \leq i \leq t$.

If $n=1$, then $G$ is a complete graph and $pc_{1,\ell}(G)=1$. Note that if $n > 1$, then $G$ is not complete so $pc_{1, \ell}(G) \geq 2$ for the remaining cases.

If $m\leq n\leq2^m$, then $K_{m,n}$ is a spanning subgraph of $G$. According to Proposition \ref{spanning} and Theorem \ref{bipartite}, we get that $pc_{1,\ell}(G)\leq pc_{1,\ell}(K_{m,n})=2$, which implies that $pc_{1,\ell}(G)=2$.

If $2\leq n<m$, we write $M_i$ for the larger one of $\sum_{j=1}^{j=i}n_j$ and $\sum_{j=i+1}^{j=t}n_{j}$, and $m_i$ for the other. We claim there must exist $i$ with $1\leq i\leq t-1$ such that $m_i\leq M_i\leq2^{m_i}$. If we suppose not, then we have $m>2^n$ since otherwise we could choose $i = t - 1$ for a contradiction. This is equivalent to $n<\frac{n}{2^n+n}(m+n)\leq\frac{1}{3}(m+n)$ since $n\geq2$. It is easily observed that there exists $2\leq i\leq t-1$ such that $M_{i-1}=\sum_{j=i}^{j=t}n_{j}$ and $M_i=\sum_{j=1}^{j=i}n_{j}$. So $M_{i-1}>2^{m_{i-1}}$ implies that $m_{i-1}<\frac{m_{i-1}}{2^{m_{i-1}}+m_{i-1}}(m+n)\leq\frac{1}{3}(m+n)$ while $M_i>2^{m_i}$ implies that $M_i>\frac{2^{m_i}}{2^{m_i}+m_i}(m+n)\geq\frac{2}{3}(m+n)$. But this leads to $n\geq n_i=M_i-m_{i-1}>\frac{1}{3}(m+n)>n$, a contradiction. So there exists $1\leq i\leq t-1$ such that $m_i\leq M_i\leq 2^{m_i}$. Therefore, $K_{m_i,M_i}$ is a spanning subgraph of $G$. Again by Proposition \ref{spanning} and Theorem \ref{bipartite}, $pc_{1,\ell}(G)=2$.

If $n>2^m$, let $U = \bigcup_{i=1}^{t-1}V_i=\{u_1,\dots,u_m\}$. We claim that two colors are not enough to color $G$ to be $(1,\ell)$-proper connected. For any $2$-coloring of $G$, by the pigeonhole principle, there exists a pair of vertices $x,y\in V_t$ such that all edges of the form $wx$ and $wy$ have the same color for all $w \in U$. Any distance $2$-proper path between $x$ and $y$ must use at least $4$ edges, meaning that it must use at least $3$ colors. We then give a $(1,\ell)$-proper-path $3$-coloring of $G$ as follows. Let $V_t=V_t^1\bigcup V_t^2$ such that $|V_t^1|=2^m$. For each element $v\in V^1_t$, we assign a vector to be $v'=(v'_1,v'_2,\dots,v'_m)$ to it such that $v'_i\in \{1,2\},~i=1,\dots,m$ and for any distinct vertices $v,w\in V^1_t$, $v'\neq w'$. Color the edge $vu_i$ with $v'_i$. For each vertex $v\in V_t^2$, we set its corresponding vector $v'=(1,2,\dots,2)$ and color the edge $vu_i$ with $v'_i$. Finally, all edges between the vertices of $\bigcup_{i=1}^{t-1}V_i$ are assigned color $3$. It is easy to check that between any pairs of vertices in $V_t^1$ or $\bigcup^{t-1}_{i=1}V_i~(i\neq j)$, there exist distance $\ell$-proper paths in between. For $x,y\in V_t^2$, clearly $xu_1zy$ is a distance $\ell$-proper $(x,y)$-path where $u_1\in V_i$, $z\in V_j$ and $1\leq i\neq j\leq t-1$. And for $x\in V_t^1$ and $y\in V_t^2$, if $x'=y'$, then it is similar to $x,y\in V_t^2$. Otherwise $xu_iy$ is a distance $\ell$-proper $(x,y)$-path, where $x'_i\neq y'_i$. Hence we have $pc_{1,\ell}(G)=3$.

Thus the proof is completed.\qed

The wheel graph $W_n$ is obtained from the cycle $C_n$ by joining a new vertex $v$ to all vertices of $C_n$ denoted by $\{u_1,\dots,u_n\}$ in the clockwise order. The vertex $v$ is the center of $W_n$. We also get results for $W_n$.

\begin{thm}\label{wheel}
Let $n\geq3$ and $\ell\geq2$, then,
\begin{eqnarray*}
pc_{1,\ell}(W_n)=\left\{
\begin{array}{rcl}
1     &      & if\ n=3,\\
2     &      & if\ 4\leq n\leq6,\\
3     &      & otherwise.
\end{array} \right.
\end{eqnarray*}
\end{thm}
\pf If $n=3$, then $W_3=K_4$, so $pc_{1,\ell}(W_3)=1$.

\begin{figure}[H]
\begin{center}
\includegraphics[scale = 0.9]{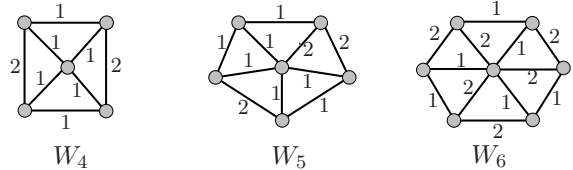}
\caption{The $(1,2)$-proper-path $2$-colorings for $W_4$, $W_5$ and $W_6$}\label{1}
\end{center}
\end{figure}

If $4\leq n\leq6$, we give $(1,2)$-proper-path $2$-colorings in Figure \ref{1}. Since $2$ colors are best possible, we know that $pc_{1,2}(W_n)=2$ for $4\leq n\leq6$. Furthermore, these $2$-colorings are also $(1,\ell)$-proper-path colorings of $W_4,~W_5,~W_6$ and so $pc_{1,\ell}(W_n)=2$ for $4\leq n\leq6$.

Now we assume that $n\geq7$. Suppose $W_n$ has a $(1,2)$-proper-path $2$-coloring $c$. Then $u_ivu_j$ must be the unique $2$-proper $(u_i,u_j)$-path if the distance between $u_i$ and $u_j$ on $C_n$ is greater than $2$. Without loss of generality, we assume $c(u_1v)=1$. Thus $c(u_iv)=2$ for $4\leq i\leq n-2$. If $n\geq9$, then $c(u_4v)=c(u_{n-2}v)=2$, a contradiction. If $n=8$, then $c(u_4v)=c(u_6v)=2$ implies that $c(u_3v)=c(u_7v)=1$, a contradiction. Finally for $n=7$, we get $c(u_{4}v) = c(u_{5}v) = 2$ so $c(u_{2}v) = c(u_{7}v) = 1$ which implies that $c(u_3v)=c(u_6v)=2$, also a contradiction. As a result, $pc_{1,\ell}(W_n)\geq pc_{1,\ell-1}(W_n)\geq\cdots\geq pc_{1,2}(W_n)>2$ for $n\geq7$. Next we give a $(1,\ell)$-proper-path $3$-coloring $f$ as follows. Set $f(u_iu_{i+1})=a\in[3]$, where $i\equiv a\ (\text{mod} \ 3)$. If $a-1\equiv b\in [3]~(\text{mod} \ 3)$, then set $f(u_iv)\in [3]\backslash\{a,b\}$, $i\in [n]\backslash \{1\}$ and set $f(u_1v)=3$. It is easy to check that $W_n(n\geq7)$ is $(1,\ell)$-proper connected in this way. So $pc_{1,\ell}(W_n)=3~(n\geq7)$, completing the proof. \qed

Additionally we also examine the $(1,\ell)$-proper connection number of the $t$-cube. The $t$-cube is the graph whose vertices are the ordered $t$-tuples of $0's$ and $1's$, two vertices being joined if and only if the $t$-tuples differ in exactly one coordinate. We denote by $Q_t$ the $t$-cube. Our result is posed below.

\begin{thm}\label{cube}
Let $t\geq1$ and $\ell\geq2$. Then
\begin{eqnarray*}
pc_{1,\ell}(Q_t)=\left\{
\begin{array}{rcl}
1~~   &      & if\ t=1,\\
2~~   &      & if\ t=2,\\
t~~   &      & if\ t\geq3\ and\ \ell\geq t,\\
\ell+1 &     & if\ t\geq3\ and\ \ell<t.
\end{array} \right.
\end{eqnarray*}
\end{thm}
\pf Since $Q_1$ is $K_2$ and $Q_2$ is $C_4$, the result can be easily verified.

Now we consider the case $t\geq3$. Since the diameter of $Q_t$ is $t$, then $pc_{1,\ell}(Q_t)\geq t$ if $\ell\geq t$. Assign the color $i$ to the edge $uv$ if $u$ and $v$ differ in the $i^{th}$ position. One can check this is a $(1,\ell)$-proper-path coloring and thus $pc_{1,\ell}(Q_t)=t$. If $\ell<t$, then $pc_{1,\ell}(Q_t)\geq\ell+1$. Next we will give a $(1,\ell)$-proper-path $(\ell+1)$-coloring, which implies $pc_{1,\ell}(Q_t)=\ell+1$. For $uv\in E(Q_t)$, let $f(uv)\in[\ell+1]$ denote the color assigned to $uv$. Let $i$ be the position where $u$ and $v$ differ and let $a \equiv i\mod \ell+1$ so that $a \in [\ell+1]$. We then define $f(uv)=a$. And the $i^{th}$ position is also called an $a$-position. Note that all edges in dimension $1$ will receive color $1$, and more generally, all edges in dimension $i$ will receive color $i\mod \ell+1$.

For any $u,v\in V(Q_t)$, assume that they differ in $p_i$ $i$-positions, $i=1,\dots,\ell+1$. Without loss of generality, we suppose that $p_{\ell+1}=\max\{p_i,i=1,\dots,\ell+1\}$. Then for arbitrary $i\in\{1,\dots,\ell+1\}$, we can find nonnegative integers $s_i$ such that $0\leq p_{\ell+1}-p_i-2s_i\leq1$. This implies that we can find $s_i$ different $i$-positions such that $u$ and $v$ share the same $i^{th}$ term. We will change each of these positions twice, making the final result unchanged in those dimensions. Then we can easily find a distance $2$-proper path from $u$ to $v$ rotating through using the colors $1,2,\dots,\ell+1$. Thus $pc_{1,\ell}(Q_t)=\ell+1$ for $t\geq3$, and the proof is completed. \qed

\section{$2$-Connected Graphs}\label{Sect:2-conn}

In this section, we mainly deal with the $2$-connected graphs and obtain an upper bound for their $(1,2)$-proper connection number.
\begin{thm}\label{2-connected}
If a graph $G$ is $2$-connected, then $pc_{1,2}(G)\leq5$.
\end{thm}
\pf Suppose $G$ is minimally $2$-connected, that is, the removal of any edge would leave $G$ not $2$-connected. We will prove this theorem by induction on the number of ears in an ear decomposition of $G$, $C_r=G_0\subset G_1\subset G_2\subset\cdots\subset G_s=G$, where $G_{i+1}=G_i\cup P^{i+1}~(0\leq i\leq s-1)$ and $P^i$ is the $i^{th}$ ear added in this ear decomposition. We use $f(e)$ to denote the color assigned to the edge $e$. Starting from a cycle, we consecutively add ears and give them appropriate colors until the graph $G$ is obtained. We denote by $start_2(P)$ the first two edges of a distance $2$-proper path $P$, and $end_2(P)$ the last two edges.

We claim that three properties hold on each stage.

1. Five colors are enough to make the graph $(1,2)$-proper connected.

\begin{figure}[H]
\begin{center}
\includegraphics[scale = 0.9]{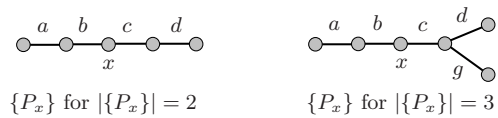}
\caption{The two forms of $\{P_x\}$}\label{2}
\end{center}
\end{figure}

2. For every vertex $x$ of the present graph, there exists a set denoted by $\{P_x\}$ consists of two or three paths of length two with a common end $x$. The two forms of $\{P_x\}$ according to its cardinalities are given respectively in Figure \ref{2}. As is shown in the figure, $\{P_x\}=\{ba,cd\}$ if $|\{P_x\}|=2$ and $\{P_x\}=\{ba,cd,cg\}$ if $|\{P_x\}|=3$. In such a way that for any $u,v$ in the present graph, there exists a distance $2$-proper $(u,v)$-path $P$ with $start_2(P)\in \{P_u\}$ and $end_2(P)\in \{P_v\}$. We refer to this path $P$ as $P_{u,v}$.

3. Let $f_{\{P_x\}}$ be the color set of edges of $\{P_x\}$. Then for each vertex $x$ in the present graph, $|f_{\{P_x\}}|\leq4$.

For the cycle $C_r$, we color its edges clockwise in the sequence of $1,2,3,1,2,\\3,\dots$. And give color $4$ to the remaining edge if $k\equiv1\ (\text{mod}\ 3)$, give colors $4$ and $5$ to the remaining two edges respectively if $k\equiv2\ (\text{mod}\ 3)$. Clearly for $\forall v\in C_k$, $v$ is the common end of two $P_3s$ on $C_r$. Then we take the set consisting of these two $P_3s$ as $\{P_v\}$. It is easy to check that the above three properties are satisfied.

Let $P=P^k$ be the $k$th ear added in the ear decomposition of $G$ and $G'=G_{k-1}$ be the graph after the removal of the internal vertices of $P$ in $G_k$. Suppose the above three properties always hold before. Because $G$ is minimally $2$-connected, $P$ has at least one internal vertex. Let $P=u(=u_1)u_2u_3\cdots u_{p+1}v(=u_{p+2})$ and $P^{-1}$ the inverse of $P$. Since $G'$ satisfies the above three properties, thus we consider $start_2(P_{u,v})=uw_1w_2\in\{P_u\}$. And we assume that the edges incident to $v$ contained in paths of $\{P_v\}$ are $vv_1$ and $vv_2$.  For any internal vertex $u_i(i\neq2,p+1)$, set $\{P_{u_i}\}=\{u_iu_{i-1}u_{i-2},u_iu_{i+1}u_{i+2}\}$. If $p=1$, set $\{P_{u_2}\}=\{u_2uw_1,u_2vv_1,u_2vv_2\}$. Otherwise, set $\{P_{u_2}\}=\{u_2uw_1,u_2u_3u_4\}$ and  $\{P_{u_{p+1}}\}=\{u_{p+1}u_pu_{p-1},u_{p+1}vv_1,u_{p+1}vv_2\}$.

Suppose that $p=1$, if $f(vv_1)=f(vv_2)$, then certainly $|f_{\{P_{u_2}\}}|\leq4$. Color the edge $u_2v$ such that $f(u_2v)\in [5]\backslash f_{\{P_v\}}$. It is possible since according to the induction, $|f_{\{P_v\}}|\leq4$. And we color the edge $u_2u$ such that $f(u_2u)\in [5]$ and $f(u_2u)\neq f(w_1w_2)$, $f(u_2u)\neq f(uw_1)$. If $f(vv_1)\neq f(vv_2)$, again we color $u_2v$ and $u_2u$ such that  $f(u_2v)\in [5] \backslash f_{\{P_v\}}$, $f(u_2u)\in\{f(e):e=u_2v,vv_1\ or\ vv_2\}$ and $f(u_2u)\neq f(w_1w_2)$, $f(u_2u)\neq f(uw_1)$. Clearly this is also possible. And in this way, we guarantee that $|f_{\{P_{u_2}\}}|\leq4$. By the induction hypothesis, for $x,y\in G'$, $P_{x,y}$ exists. For $u_2,v$, $u_2uP_{u,v}v$ are as required, and for $u_2$ and $w\neq v$, $u_2vP_{v,w}w$ is the path we demand.

Suppose that $p\geq2$, if $f(vv_1)=f(vv_2)$, again $|f_{
\{P_{u_{p+1}}\}}|\leq4$. Color the edges $u_{p+1}v$ and $u_pu_{p+1}$ such that $f(u_{p+1}v)\in [5]\backslash f_{\{P_v\}}$ and $f(u_pu_{p+1})\notin\{f(e):e=u_{p+1}v,\ vv_1\ or\ vv_{2}\}$. As for any uncolored edge on $P$, it is assigned some color different from its two neighboring edges from each side(four edges altogether) on the path $w_2w_1uPv$. This is possible since we have five distinct colors. If $f(vv_1)\neq f(vv_2)$, again we color $u_{p+1}v$ such that $f(u_{p+1}v)\in [5]\backslash f_{\{P_v\}}$. When $p=2$, we set $f(uu_2)\neq f(uw_1)$, $f(uu_2)\neq f(w_1w_2)$ and $f(uu_2)\in \{f(e):e=u_3v,vv_1~\text{or}~vv_2\}$. For other conditions, we set $f(u_{p-1}u_p)\neq f(uw_1)$ and $f(u_{p-1}u_p)\in\{f(e):e=vv_1\ or\ vv_{2}\}$. Then we color the edge $u_pu_{p+1}$ such that $f(u_pu_{p+1})\in [5]\backslash \{f(e):e=u_{p+1}v,vv_1,vv_2,u_{p-1}u_p,u_{p-2}u_{p-1}\}$. And for any other uncolored edge on $P$, we give it a color different from its two neighboring edges from each side on the path $w_2w_1uPv$. So the coloring is done. For $x,y\in G'$, $P_{x,y}$ exists by induction. For $u_i,v~(2\leq i\leq p+1)$, the path $u_iP^{-1}uP_{u,v}$ is as required, and for $u_i~(2\leq i\leq p+1)$ and $w\in G'\backslash v$, $u_iPvP_{v,w}$ is the path we demand. For $u_i,u_j~(2\leq i<j\leq p+1)$, their distance $2$-proper path is $u_iP^{-1}uP_{u,v}vu_{p+1}P^{-1}u_j$. It can be easily verified that the three properties hold.

Therefore, for any minimally $2$-connected graph $G$, we have $pc_{1,2}(G)\leq5$. Finally we consider the situation that $G$ is not minimally $2$-connected. We choose a minimally $2$-connected subgraph $H\subset G$. Together with the above analysis and Proposition \ref{spanning}, we obtain that $pc_{1,2}(G)\leq5$. And thus the proof is completed.\qed

\section{Graph Operations}\label{Sect:Oper}
In this section, we consider the $(1,2)$-proper connection number of graphs obtained by some graph operations. These operations include the joins of graphs, cartesian products of graphs and permutation graphs.

The \emph{join} $G\vee H$ of two graphs $G$ and $H$ has vertex set $V(G)\cup V(H)$ and edge set $E(G)\cup E(H)\cup\{uv:u\in V(G)~\text{and}~v\in V(H)\}$.

\begin{thm}\label{join}
If $G$ and $H$ are nontrivial connected graphs, then $pc_{1,2}(G\vee H)\leq 3$.
\end{thm}
\pf Since $G$ must contain a spanning complete bipartite graph $K_{m,n}$ with $m,n\geq2$. Then according to Theorem \ref{bipartite} and Proposition \ref{spanning}, the result is clear.\qed

The \emph{Cartesian product} $G\square H$ of two graphs $G$ and $H$ has vertex set $V(G\square H)=V(G)\times V(H)$ and two distinct vertices $(u,v)$ and $(x,y)$ of $G\square H$ are adjacent if either $ux\in E(G)$ and $v=y$ or $vy\in E(H)$ and $u=x$. We show next that the Cartesian product of graphs is a special graph class with rather small $(1,2)$-proper connection numbers.


\begin{thm}\label{cartesian}
Let $G$ and $H$ be nontrivial connected graphs, and not both of them are complete graphs. Then

\noindent$(\romannumeral 1)$ $pc_{1,2}(G\square H)=3$ except if $G~(\text{or}~H)$ is a star and every spanning tree of $H~(\text{respectively}~G)$ has radius at least $3$.

\noindent$(\romannumeral 2)$ $pc_{1,2}(G\square H)\leq 4$ if $G~(\text{or}~H)$ is a star and every spanning tree of $H~(\text{respectively}~G)$ has radius at least $3$.
\end{thm}

\pf Let $S$ and $T$ be spanning trees of $G$ and $H$, respectively. Let $V(S)=\{u_1,u_2,\dots,u_m\}$ and $V(T)=\{v_1,v_2,\dots,v_n\}$. And we take $u_1$ and $v_1$ as the roots of $S$ and $T$ respectively. In the Cartesian product $G\square H$, we denote by $S_i$ the tree $S$ corresponding to the vertex $v_i$ of $T$, $S_{i,j}$ the vertex $u_j$ for $S_i$. And $T_i$, $T_{i,j}$ are defined similarly. We also write $P_{S_{i,j}}$ for the unique path in the tree $S_i$ from the root $S_{i,1}$ to the vertex $S_{i,j}$, and $P_{S_{i,j}}^{-1}$ its reverse. And $P_{T_{s,t}}$, $P_{T_{s,t}}^{-1}$ have similar definitions. Also write $\ell(P)$ for the length of the path $P$.

Notice that since not both $G$ and $H$ are complete graphs, then we have $diam(G\square H)\geq3$. Thus $pc_{1,2}(G\square H)\geq3$.

\emph{Proof for $(\romannumeral 1)$}: Suppose first that neither $G$ nor $H$ is $K_3$. By the condition of $(\romannumeral 1)$, we can always choose $S,~T,~u_1$ and $v_1$ in such a way that one of the following holds
\begin{description}
\item{(a)\hspace{0.17cm}} $ecc_S(u_1)=2$ (or $ecc_T(v_1)=2$) and \\ $ecc_T(v_1)\leq2$ (respectively $ecc_S(u_1)\leq2$), or
\item{(b)\hspace{0.17cm}} $ecc_S(u_1)\geq3$, and $ecc_T(v_1)\geq3$.
\end{description}
Since $S\square T$ is a spanning subgraph of $G\square H$, it suffices to provide a $(1,2)$-proper-path $3$-coloring for $S\square T$. First we color the edges of $T_1$ such that for all $t$ with $1\leq t\leq n$, the path $P_{T_{1,t}}$ in $T_1$ is a distance $2$-proper path. Then we assign colors to the edges of $S_1$ so that for $j$ with $1\leq j\leq m$ and for $t$ with $1\leq t\leq n$, the path $P_{S_{1,j}}^{-1}P_{T_{1,t}}$ is a distance $2$-proper path. For $2\leq i\leq n$, we give colors to the edges of $S_i$ such that for all $j_{1}$ with $1\leq j_1\leq m$ and all $j_{i}$ with $1\leq j_i\leq m$, the path $P_{S_{1,j_1}}^{-1}P_{T_{1,i}}P_{S_{i,j_i}}$ is $2$-proper. Finally we color the edges of $T_s$ (for $2\leq s\leq m$) in such a way that for $1\leq t_1\leq n$ and $1\leq t_s\leq n$, the path $P_{T_{1,t_1}}^{-1}P_{S_{1,s}}P_{T_{s,t_s}}$ is $2$-proper.

We then verify that the coloring given above is a $(1,2)$-proper-path coloring. We distinguish four cases to analyze based on the locations of the two vertices that we are trying to connect:

1. Consider vertices $S_{1,i}$ and $S_{1,j}$ (or similarly $T_{1,i}$ and $T_{1,j}$) with $i\neq j$. If one vertex, say $S_{1, i}$, lies on the unique path ($P_{S_{1,j}}$) from the other ($S_{1, j}$) to the root, the path between them within $S_{1}$ is trivially the desired distance $2$-proper path. Otherwise, without loss of generality, we suppose $j > 1$ so regardless which of (a) or (b) holds, we can find a vertex $T_{1,k}$ such that $\ell(P_{S_{1,j}})+\ell(P_{T_{1,k}})=0\ (\text{mod}\ 3)$. Since $P_{T_{j,k}}^{-1}P_{S_{1,j}}^{-1}P_{T_{1,k}}P_{S_{k,j}}$ is a distance $2$-proper path and it has length $\ell=0\ (\text{mod}\ 3)$, this means that the path $P_{T_{1,k}}P_{S_{k,j}}P_{T_{j,k}}^{-1}$ is also a distance $2$-proper path. Therefore, the path $P_{S_{1,i}}^{-1}P_{T_{1,k}}P_{S_{k,j}}P_{T_{j,k}}^{-1}$ is distance $2$-proper, and the one we desire. The case for $T_{1,i}$ and $T_{1,j}$ can be dealt with similarly.

2. Consider vertices $S_{s,i}$ and $S_{s,j}$ for $(2\leq s\leq n)$ (or similarly $T_{t,i}$ and $T_{t,j}$ for $(2\leq t\leq m)$). For the pair $S_{s,i}$ and $S_{s,j}$, if one vertex lies in the path from the root to the other one, then the unique path between them in $S_{s}$ is the desired distance $2$-proper path. Otherwise the $(S_{s,i},S_{s,j})$-path $P_{S_{s,i}}^{-1}P_{T_{1,s}}^{-1}P_{S_{1,j}}P_{T_{j,s}}$ is clearly distance $2$-proper. With the same method, we can easily check the case for $T_{t,i}$ and $T_{t,j}$.

3. Consider vertices $S_{1,i}$ and $S_{s,j}$ for $(1\neq s,i\neq j)$ (or similarly $T_{1,i}$ and $T_{s,j}$ for $(1\neq i,i\neq j)$). In between $S_{1,i}$ and $S_{s,j}$, the path $P_{S_{s,j}}^{-1}P_{T_{1,s}}^{-1}P_{S_{1,i}}$ is the desired distance $2$-proper path between them. The case for $T_{1,i}$ and $T_{s,j}$ can also be checked easily.

4. Consider vertices $S_{i,j}$ and $S_{s,t}$ for $(i,j,s,t\neq 1,i\neq s,j\neq t)$. For the vertices $u_j$ and $u_t$ in $S$, if one lies in the path from the root to the other, then without loss of generality, we let $u_j$ be the one closest to the root. In this way, the path $P_{S_{i,j}}^{-1}P_{T_{1,i}}^{-1}P_{S_{1,t}}P_{T_{t,s}}$ is the desired distance $2$-proper path between $S_{i,j}$ and $S_{s,t}$.

Next we may assume that one of $G$ and $H$ is $K_3$. Without loss of generality, we assume that $H$ is $K_3$ and $V(H)=\{t_1,t_2,t_3\}$. Then similarly to above, let $S$ be a spanning tree of $G$ and $V(S)=\{u_1,u_2,\dots,u_m\}$. In the Cartesian product $G\square H=G\square K_3$, we denote by $S_i$ the tree $S$ corresponding to the vertex $t_i$ of $H$ and $S_{i,j}$ the vertex $u_j$ of $S_i$. And $H_i$, $H_{i,j}$ are defined similarly. We also write $P_{S_{i,j}}$ for the unique path in the tree $S_i$ from the root $S_{i,1}$ to the vertex $S_{i,j}$, and $P_{S_{i,j}}^{-1}$ its reverse.

Similar to the above analysis, we only need to give a $(1,2)$-proper-path $3$-coloring for $S\square H$. First we color the edges of $S_2$ such that for $1\leq i\leq m$, the path $P_{S_{2,i}}$ in $S_2$ is a distance $2$-proper path. Then we assign colors to the edges of $S_1$ and $S_3$ such that for all $j$ and $k$ with $1\leq j\leq m$ and $1\leq k\leq m$, the path $P_{S_{i,j}}^{-1}S_{i,1}S_{2,1}P_{S_{2,k}}~(i=1,3)$ is a distance $2$-proper path. For the edge $S_{1,i}S_{2,i}~(2\leq i\leq m)$, we give it an appropriate color such that $P_{S_{2,i}}S_{1,i}~(2\leq i\leq m)$ is a distance $2$-proper path. And for the edge $S_{2,i}S_{3,i}~(2\leq i\leq m)$, we give it a suitable color such that $P_{S_{3,i}}S_{2,i}~(2\leq i\leq m)$ is a distance $2$-proper path. As for the edge $H_{1,1}H_{1,3}$, we give it an appropriate color such that $P_{S_{1,m}}^{-1}H_{1,3}$ is a distance $2$-proper path. For the edge $S_{1,i}S_{3,i}~(2\leq i\leq m)$, we give it a suitable color such that $S_{3,1}S_{2,1}P_{S_{2,i}}S_{1,i}S_{3,i}$ is a distance $2$-proper path. It can be verified that for any pair of vertices, there exists a distance $2$-proper path between them and we omit the details here.

\emph{Proof for $(\romannumeral 2)$}: Assume that $G$ is a star and all spanning trees of $H$ have radius at least $3$. Then $S=G$ is a star and $T$ is a spanning tree of $H$ with vertex set $\{u_1,u_2,\dots,u_m\}$ and $\{v_1,v_2,\dots,v_n\}$, respectively. Take the center vertex $u_1$ as the root of $S$ and one end vertex $v_1$ of a longest path of $T$ as its root. The symbols $S_{i,j},T_{i,j},P_{S_{i,j}},P_{T_{s,t}},P^{-1}_{S_{i,j}}$ and $P^{-1}_{T_{s,t}}$ are all defined as the beginning of the proof.

\begin{figure}[H]
\begin{center}
\includegraphics[scale = 0.9]{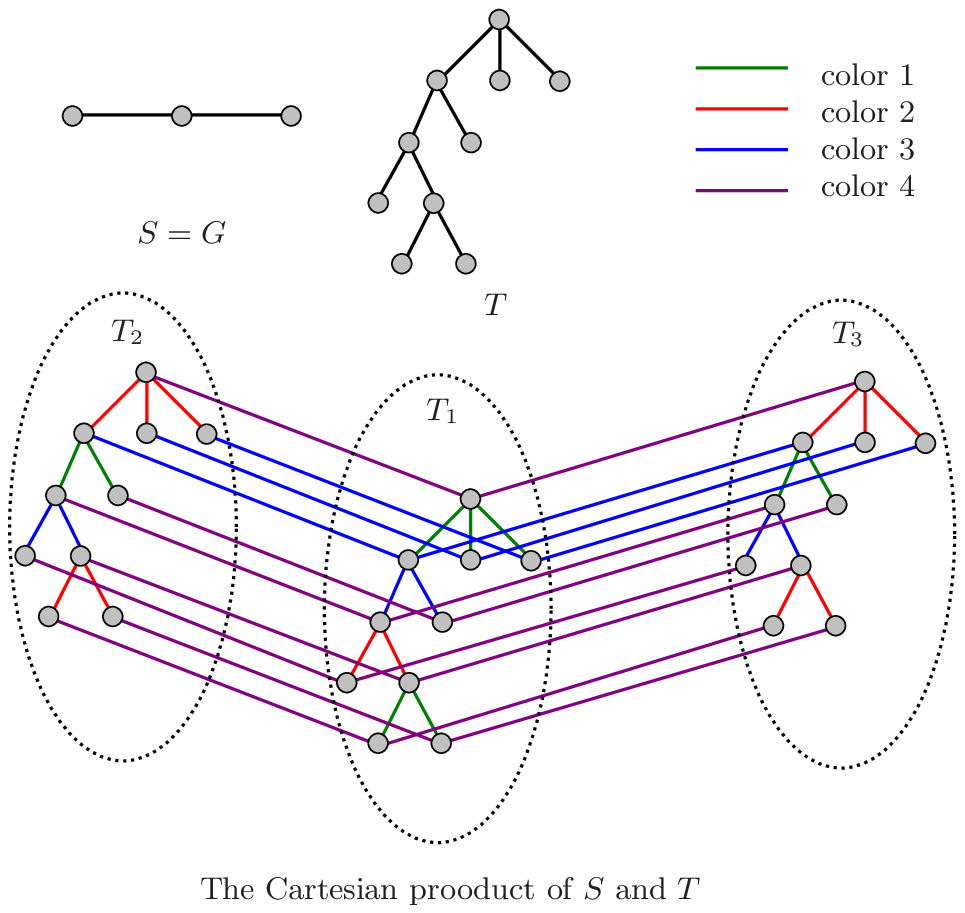}
\caption{An example for Theorem \ref{cartesian}$(\romannumeral 2)$}\label{3}
\end{center}
\end{figure}

Now we give a $(1,2)$-proper-path $4$-coloring for $S\square T$, which implies the conclusion in $(\romannumeral 2)$. 
First we give colors $1,2,3$ to the edges of $T_1$ such that for all $j$ with $1\leq j\leq n$, the path $P_{T_{1,j}}$ is a distance $2$-proper path. Then assign the color $4$ to all edges of $S_1$ and use $1,2,3$ to color the edges of $T_2$ such that for all $j$ and $t$ with $1\leq j\leq n$ and $1\leq t\leq n$, we have $P_{T_{1,j}}^{-1}T_{1,1}T_{2,1}P_{T_{2,t}}$ is a distance $2$-proper path. Then give the edges of $T_i~(3\leq i\leq m)$ the same colors as the corresponding edges in $T_2$. For $2\leq r\leq n$, we color the edges of $S_r$ such that for all $i$ with $2\leq i\leq m$, the cycles $P_{T_{1,r}}^{-1}P_{S_{1,i}}P_{T_{i,r}}P_{S_{r,i}}^{-1}$ and $P_{S_{r,i}}^{-1}P_{T_{1,r}}^{-1}P_{S_{1,i}}P_{T_{i,r}}$ are distance $2$-proper. An example of such a coloring is depicted in Figure~\ref{3}. Note that for each $r$ with $1\leq r\leq n$, all of the edges in $S_r$ share the common color. This can be verified as a $(1,2)$-proper-path $4$-coloring of $S\square T$ but we omit the details.\qed

\begin{remark}
The second part of the above theorem shows that $4$ is an upper bound of $pc_{1,2}(G\square H)$ for the case $G~(\text{or}~H)$ is a star and every spanning tree of $H~(\text{or}~G)$ has radius at least $3$. However, it is not clear whether this upper bound is sharp or it can be decreased to $3$.
\end{remark}

We also consider the $(1,\ell)$-proper connection number of permutation graphs. Let $G$ be a graph with $V(G)=\{v_1,\dots,v_n\}$ and $\alpha$ be a permutation of $[n]$. Let $G'$ be a copy of $G$ with vertices labeled $\{u_{1}, \dots, u_{n}\}$ where $u_{i} \in G'$ corresponds to $v_{i} \in G$. Then the \emph{permutation graph} $P_\alpha(G)$ of $G$ can be obtained from $G \cup G'$ by adding all edges of the form $v_i u_{\alpha(i)}$. Next we present our result on permutation graphs of traceable graphs.
\begin{thm}\label{permutation graph}
Let $G$ be a nontrivial traceable graph of order $n$, then
\begin{eqnarray*}
pc_{1,\ell}(P_\alpha(G))\leq\ell+1
\end{eqnarray*}
for each permutation $\alpha$ of $[n]$.
\end{thm}
\pf We use $f(e)~(e\in E(G))$ to represent the color assigned to $e$. Let $P=v_1v_2\cdots v_n$ be a hamiltonian path of $G$. Then $P'=u_1u_2\cdots u_n$ is a hamiltonian path of $G'$. Besides, we write $P^{-1}$ and $P'^{-1}$ the reverse of $P$ and $P'$, respectively. Firstly we consider the cases when $\alpha(n)=1~or~n$, then clearly $P_\alpha(G)$ is traceable and the theorem holds according to Proposition \ref{traceable}.

Otherwise, we suppose $\alpha(n)=i~(2\leq i\leq n-1)$. We color the edges of $P$ with $\ell+1$ colors following the sequence $1,2,\dots,\ell+1,1,2,\dots,\ell+1,\dots$. We then color the remaining edges in the three paths $v_1Pv_nu_iP'^{-1}u_1$,  $v_1Pv_nu_iP'u_n$ and $u_{\alpha(1)}v_1Pv_n$ so that each follows the sequence $1,2,\dots,\ell+1,1,2,\dots,\ell+1,\dots$. Finally set $f(v_ju_{\alpha(j)})=f(v_{j-1}v_j)~(2\leq j\leq n-1)$. In this way, it is easy to see the distance $\ell$-proper paths between all pairs of vertices except between $u_s$ and $u_t$ with $1\leq s\leq i-1$ and $i+1\leq t\leq n$. In this case, the path $u_sP'u_iv_nP^{-1}v_{\alpha^{-1}(t)}u_t$ is the desired distance $\ell$-proper path. Thus the proof is complete.\qed


\begin{thebibliography}{1}

\bibitem{ALLZ}
E. Andrews, E. Laforge, C. Lumduanhom and P. Zhang, On proper-path colorings in graphs, {\it J. Combin. Math. Combin. Comput.} to appear.

\bibitem{BFG}
V. Borozan, S. Fujita, A. Gerek, C. Magnant, Y. Manoussakis, L. Montero,  Zs. Tuza, Proper connection of graphs, {\it Discrete Math.} {\bf 312(17)}(2012), 2550--2560.

\bibitem{BM}
J. A. Bondy, U. S. R. Murty, Graph Therory, GTM 244, Springer-Verlag, New York, 2008.

\bibitem{CJMZ}
G. Chartrand, G. L. Johns, K. A. McKeon, P. Zhang, Rainbow connection in graphs, {\it Math Bohemica.} {\bf 133(1)}(2008), 5--98.

\bibitem{FJ}
J. L. Fouquet, J. L. Jolivet, Strong edge-coloring of graphs and applications to multi-$k$-gons, {\it Ars Combin.} {\bf 16A}(1983), 141--150.

\bibitem{KY}
M. Kriveleich, R. Yuster, The rainbow connection of a graph is (at most) reciprocal to its minimum degree, {\it J. Graph Theory} {\bf 71}(2012), 206--218.

\bibitem{LM}
X. Li, C. Magnant, Properly colored notations of connectivity-a dynamic survey, Theory and Applications of Graphs, {\bf 0(1)}(2015), Article 2.

\bibitem{LS}
X. Li, Y. Sun, Rainbow Connections of Graphs, Springer Briefs in Math., Springer, New York, 2012.

\bibitem{LS2}
X. Li, Y. Shi, Rainbow connection in 3-connected graphs, {\it Graphs} \& {\it Combin.} {\bf 29(5)}(2013), 1471--1475.

\bibitem{LSS}
X. Li, Y. Shi, Y. Sun, Rainbow connections of graphs: A survey, {\it Graphs} \& {\it Combin.} {\bf 29}(2013), 1--38.
\end{thebibliography}
\end{document}